\documentclass[12pt,twoside,a4paper]{amsart}
\usepackage[utf8]{inputenc}
\usepackage{xspace}
\usepackage{latexsym}
\usepackage{amsfonts}
\usepackage{mathrsfs}
\usepackage{amssymb}
\usepackage{amsmath}
\usepackage{amsxtra}
\usepackage{amsthm}
\usepackage{amsaddr}
\usepackage[geometry]{ifsym}
\usepackage{graphics}
\usepackage{titlesec}
\usepackage[french,english]{babel}
\usepackage[bookmarks]{hyperref}
\usepackage{enumerate}
\usepackage{color}

\usepackage{lineno}\nolinenumbers

\titleformat{\section}[block]
   {\filcenter\normalfont}
   {\S \thesection.}{.5em}{\rm\protect\uppercase}
\titlespacing{\section}
   {0pc}{*4}{*1}

\titleformat{\subsection}[block]
   {\bf}
   {\bf\thesubsection.}{0.5em}{}
\titlespacing{\subsection}
   {\parindent}{1.5ex plus .1ex minus .2ex}{.5em}

\titleformat{\subsubsection}[block]
   {\normalfont\it}
   {\thesubsubsection.}{.5em}{}
\titlespacing{\subsubsection}
   {\parindent}{1.5ex plus .1ex minus .2ex}{.5em}

\renewcommand{\thesubsection}{\arabic{subsection}}

\newcounter{aninhamento}
\renewcommand{\theaninhamento}{\alph{aninhamento}}
\newcounter{countlistaa}
\newcounter{countlistab}[countlistaa]
\newcounter{countlistac}[countlistab]
\newcounter{countlistad}[countlistac]
\newcounter{countlistae}[countlistad]
\newcounter{countlistaf}[countlistae]

{\addtocounter{aninhamento}{1}%
 \renewcommand{\thecountlistaa}{#1{countlistaa}}%
 \renewcommand{\thecountlistab}{#1{countlistab}}%
 \renewcommand{\thecountlistac}{#1{countlistac}}%
 \renewcommand{\thecountlistad}{#1{countlistad}}%
 \renewcommand{\thecountlistae}{#1{countlistae}}%
 \renewcommand{\thecountlistaf}{#1{countlistaf}}%
 \begin{list}{\rm(#1{countlista\theaninhamento})}{\usecounter{countlista\theaninhamento} \setlength{\leftmargin}{0pt}}}%
{\end{list}\addtocounter{aninhamento}{-1}}

{\addtocounter{aninhamento}{1}%
 \renewcommand{\thecountlistaa}{#1{countlistaa}}%
 \renewcommand{\thecountlistab}{#1{countlistab}}%
 \renewcommand{\thecountlistac}{#1{countlistac}}%
 \renewcommand{\thecountlistad}{#1{countlistad}}%
 \renewcommand{\thecountlistae}{#1{countlistae}}%
 \renewcommand{\thecountlistaf}{#1{countlistaf}}%
 \begin{list}{\rm(#1{countlista\theaninhamento})}{\usecounter{countlista\theaninhamento}}}%
{\end{list}\addtocounter{aninhamento}{-1}}

{\addtocounter{aninhamento}{1}%
 \renewcommand{\thecountlistaa}{#1{countlistaa}}%
 \renewcommand{\thecountlistab}{#1{countlistab}}%
 \renewcommand{\thecountlistac}{#1{countlistac}}%
 \renewcommand{\thecountlistad}{#1{countlistad}}%
 \renewcommand{\thecountlistae}{#1{countlistae}}%
 \renewcommand{\thecountlistaf}{#1{countlistaf}}%
 \begin{list}{(#2#1{countlista\theaninhamento})}{\usecounter{countlista\theaninhamento}\setlength{\leftmargin}{0pt}}}%
{\end{list}\addtocounter{aninhamento}{-1}}

{\addtocounter{aninhamento}{1}%
 \renewcommand{\thecountlistaa}{#1{countlistaa}}%
 \renewcommand{\thecountlistab}{#1{countlistab}}%
 \renewcommand{\thecountlistac}{#1{countlistac}}%
 \renewcommand{\thecountlistad}{#1{countlistad}}%
 \renewcommand{\thecountlistae}{#1{countlistae}}%
 \renewcommand{\thecountlistaf}{#1{countlistaf}}%
 \begin{list}{(#2#1{countlista\theaninhamento})}{\usecounter{countlista\theaninhamento}}}%
{\end{list}\addtocounter{aninhamento}{-1}}

\newtheorem{theorem}{\rm\scshape{Theorem}}
\newtheorem{corollary}{\rm\scshape{Corollary}}
\newtheorem{rem}{\it Remark}
\newtheorem{expl}{\it Example}

\newcommand{\M}{\vdd{M}}

\newcommand{\Hor}{\ensuremath{\mathrm{Hor}}}

\newcommand{\Ver}{\ensuremath{\mathrm{Ver}}}

\newcommand{\funcAB}[3]{\ensuremath{#1:#2\rightarrow #3}}

\newcommand{\lb}[2]{\ensuremath{[#1,#2]}}                                                % Lie bracket
\newcommand{\lins}[1]{[#1]\xspace}                                                       % linear span
\newcommand{\pullb}[1]{\ensuremath{{#1}^*}}                                              % pull back
\newcommand{\pushf}[1]{\ensuremath{{#1}_*}}                                              % push forward

\newcommand{\connection}[2]{\ensuremath{\nabla_{\textstyle #1}#2}}
\newcommand{\indconnection}[3]{\ensuremath{\nabla_{\textstyle #1}^{#3}#2}}
\newcommand{\connector}{\ensuremath{\kappa}\xspace}
\newcommand{\Sg}{\tens{S}}                        %geodesic spray
\newcommand{\lv}{\ensuremath{\lambda}\xspace}     %vertical lift
\newcommand{\lh}{\ensuremath{\mathrm{H}}\xspace}  %horizontal lift
\newcommand{\pardev}{\ensuremath{\mathbb{P}}\xspace}

\newcommand{\fibdev}{\ensuremath{\mathbb{F}}\xspace}

\newcommand{\Nn}{\ensuremath{\mathbb{N}}\xspace}

\newcommand{\Rn}{\ensuremath{\mathbb{R}}\xspace}

\newcommand{\gtensor}{\ensuremath{\tens{g}}\xspace}                             % Metric Tensor
\newcommand{\cv}{\ensuremath{\mathfrak{X}}}                                   % set of vector fields in a manifold
\newcommand{\ft}{\ensuremath{\tens{T}}}                                         % tangent bundle
\newcommand{\et}[1]{\ensuremath{\tens{T}_{#1}}}                                 % tangent space
\newcommand{\secao}[1][\infty]{\ensuremath{\Gamma^{#1}}}                        % sections of a bundle
\newcommand{\secloc}[1][\infty]{\ensuremath{\Gamma^{#1}_{\mathsf{loc}}}}                        % local sections of a bundle

\newcommand{\D}{\ensuremath{\tens{T}}}                                          % tangent map
\newcommand{\Ck}[1]{\ensuremath{\tens{C^{#1}}}\xspace}                          % functions of class C #1
\newcommand{\vf}{\ensuremath{\vec{X}}}                                            % Hamiltonian vector field
\newcommand{\lin}{\ensuremath{\mathsf{L}}\xspace}                               % linear maps
\DeclareMathOperator{\rk}{rk\,}                                                 % rank
\DeclareMathOperator{\dom}{dom\,}                                               % domain
\DeclareMathOperator{\linspan}{span\,}                                          % span
\DeclareMathOperator*{\sd}{\oplus}                                              % direct sum

\DeclareMathOperator{\shfg}{\mathcal{C}}                                        % sheaf of germs of sections

\newcommand{\ct}{\ensuremath{\mathsf{R}}\xspace}                                % curvature tensor

\newcommand{\vincl}{\ensuremath{\mathscr{D}}\xspace}                            % linear constraint
\newcommand{\Fcal}{\ensuremath{\mathcal{F}}}

\newcommand{\tens}[1]{\ensuremath{\mathsf{#1}}\xspace}
\renewcommand{\vec}[1]{\ensuremath{\mathrm{#1}}\xspace}
\newcommand{\vdd}[1]{\ensuremath{\tens{#1}}\xspace}                           % Manifold
\newcommand{\setcal}[1]{\ensuremath{\mathcal{#1}}\xspace}                       % set
\newcommand{\shf}[1]{\ensuremath{\mathcal{#1}}\xspace}                            % sheaf 

\newcommand{\N}{\vdd{N}}

\newcommand{\setU}{\setcal{U}}

\newcommand{\setV}{\setcal{V}}

\newcommand{\Dshf}{\shf{D}}
\newcommand{\Gsf}{\mathsf{G}}
\newcommand{\Psf}{\mathsf{P}}

\newcommand{\bigpar}[1]{\ensuremath{\bigl(}#1\ensuremath{\bigr)}}

\newcommand{\talque}{\ensuremath{\mid}}

\title[Smoothing Problem in Chow's Theorem]{A Note on the Smoothing Problem in Chow's Theorem}

\author{Waldyr M. Oliva}
\address{Universidade de São Paulo\\Instituto de Matemática e Estatística\\Departamento de Matemática Aplicada\\Rua do Matão, 1010, 05508-090 São Paulo, Brazil.}
\address{Universidade de Lisboa\\Instituto Superior T\'ecnico\\ Center for Mathematical Analysis, Geometry and Dynamical Systems, Departamento de Matemática, Av. Rovisco Pais, 1049-001 Lisbon, Portugal.}
\email{wamoliva@math.ist.utl.pt}

\author{Gláucio Terra}
\address{Universidade de São Paulo\\Instituto de Matemática e Estatística\\Departamento de Matemática\\Rua do Matão, 1010, 05508-090 São Paulo, Brazil.}
\email{glaucio@ime.usp.br (corresponding author)}

\begin{document}

\maketitle

\bigskip
\begin{center}
	The authors offer the present paper to Prof. Giorgio Fusco for many years of collaboration and friendship.
\end{center}

\begin{abstract}
This paper concerns a solution of the smoothing problem in Chow-Rashevskii's connectivity theorem proposed in \cite{BryantHsu}.
\end{abstract}

\section{Introduction and Objectives}

Let $\M$ be a finite dimensional paracompact smooth manifold endowed with a smooth linear subbundle $\vincl$ of $\ft\M$. The well-known 
Chow-Rashevskii's connectivity theorem (see \cite{Chow} and generalizations by P. Stefan in \cite{Stefan}, \cite{Stefan2} and by H. Sussmann in \cite{Sussmann}) asserts that, if $\vincl$ is bracket-generating, any two points in the same connected component of $\M$ may be connected by a sectionally smooth path tangent to $\vincl$. The question of whether or not any two points in $\M$ may be connected by a smooth horizontal immersion was posed by R. Bryant and L. Hsu in \cite{BryantHsu} and affirmatively answered by M. Gromov in \cite{Gromov}, who named the problem as ``the smoothing problem in Chow's theorem''.

The purpose of this note is to present an alternative approach to Gromov's solution by means of a method that, to our taste, seems to be more geometrically intuitive. Besides, it conveys some additional information on the connectivity problem: we prove in theorem \ref{thm:ChowSmooth} and its corollary \ref{cor:ChowSmooth} that, if the distribution $\vincl$ is bracket-generating, any two points in a connected open set $\setU\subset\M$ may be connected on $\setU$ by a smooth horizontal $\Ck{1}$-immersion with arbitrary given initial and final velocities in $\vincl$. Our method is quite simple: given $p,q\in\setU$, $v_p\in\vincl_p\setminus\{0\}$ and $v_q\in\vincl_q\setminus\{0\}$, we apply the orbit theorem to show that $v_p$ and $v_q$ may be connected on $(\vincl\rvert_\setU)^\ast$ (i.e. $\vincl\rvert_\setU$ with the zero section removed) by means of a sectionally smooth curve whose smooth arcs are integral curves of second order vector fields on $\vincl$, i.e. local smooth sections of $\funcAB{\tau_\vincl}{\ft\vincl}{\vincl}$ whose integral curves are lifts of smooth curves on $\M$. It then follows that the projection on $\M$ of this sectionally smooth curve is a horizontal $\Ck{1}$ immersed curve connecting $p$ and $q$ on $\setU$, whose initial and final velocities coincide with $v_p$ and $v_q$, respectively. This method may also be applied in case the linear subbundle $\vincl$ is not bracket-generating: we prove in theorem \ref{thm:SussmannSmooth} that, if $\vincl$ satisfies Sussmann's necessary and sufficient condition for reachability given in theorem 7.1 of \cite{Sussmann}, then any two points in the same connected component of $\M$ may be connected by smooth horizontal $\Ck{1}$-immersion with arbitrary given initial and final velocities in $\vincl$.

\section{Preliminaries and Notation}

\subsection{Smooth Distributions}
We denote the tangent bundle of a finite dimensional paracompact smooth\footnote{\emph{smooth} in this paper means ``$\Ck{\infty}$''} manifold $\M$ by $\funcAB{\tau_\M}{\ft\M}{\M}$. 
Following the notation and definitions in \cite{Sussmann}, a \emph{distribution} $\vincl$ on $\M$ is a family $\{\vincl_x\}_{x\in\M}$ of linear subspaces of each fiber of the tangent bundle $\funcAB{\tau_\M}{\ft\M}{\M}$. The distribution $\vincl$ is called \emph{smooth} if $\vincl_x$ varies smoothly with $x\in\M$, in the sense that there exists a set $\Dshf$ of locally defined smooth vector fields on $\M$ such that, for each $x\in\M$, $\vincl_x = \linspan\{V(x)\talque V\in\Dshf, x\in\dom V\}$ (with the convention that the linear span of the empty set is $\{0\}$). If that is the case, we say that the smooth distribution $\vincl$ is generated by $\Dshf$. Equivalently, and perhaps more naturally, the distribution $\vincl$ is smooth if there exists a subsheaf $\Dshf$ of the sheaf $\shfg^\infty_{\ft\M}$ of germs of smooth sections of $\ft\M$ (considered as a sheaf of $\Ck{\infty}(\M)$-modules) such that, for each $x\in\M$, $\vincl_x=\{V(x)\talque V\in\Dshf_x\}$ (where $\Dshf_x$ denotes the stalk of $\Dshf$ over $x$). We avoid, however, the use of sheaves, in order to keep the notation and formalism compatible with that of \cite{Sussmann} and \cite{Stefan}, \cite{Stefan2}. Note that the rank of $\vincl_x$ depends on $x$, i.e. $\vincl$ need not be a linear subbundle of $\ft\M$ (but we do assume that as a hypothesis for our main results). If $\Dshf$ is a set of locally defined smooth vector fields on $\M$, we denote by $\lins{\Dshf}$ the smooth distribution generated by $\Dshf$.

We say that $V$ is a (local) smooth section of a smooth distribution $\vincl$ if it is a smooth (local) section of $\funcAB{\tau_\M}{\ft\M}{\M}$ defined on an open set $\setU\subset\M$ such that, for all $x\in\setU$, $V(x)\in\vincl_x$. We denote the set of such local smooth sections by $\secloc(\vincl)$; it is clear that the smooth distribution $\vincl$ is generated by $\secloc(\vincl)$. 

Given two locally defined smooth vector fields on $\M$, their Lie bracket is a well-defined smooth vector field on the intersection of their domains. We say that a set of locally defined smooth vector fields $\Dshf$ on $\M$ is \emph{involutive} if it is closed by the operation of taking Lie brackets. Any set $\Dshf$ of locally defined smooth vector fields $\Dshf$ on $\M$ is contained in a smallest involutive set of locally defined smooth vector fields on $\M$, which we denote by $\Dshf_\ast$. Indeed, the family $\Fcal$ of all involutive sets of locally defined smooth vector fields containing $\Dshf$ is nonempty (since $\secloc(\ft\M)$ is such a set) and $\cap\Fcal$ does the work. We say that a smooth distribution $\vincl$ on $\M$ is \emph{involutive} if so is $\secloc(\vincl)$.

We say that a smooth distribution $\vincl$ on $\M$ is \emph{bracket-generating}
if the smooth distribution generated by $\secloc(\vincl)_\ast$ coincides with $\ft\M$.

\subsection{Orbits of Local Groups of Diffeomorphisms and Distributions}

A \emph{local group of diffeomorphisms} $\Gsf$ on $\M$ is a set of smooth diffeomorphisms defined on open subsets of $\M$ which is closed under compositions and under taking inverses, i.e. if $\funcAB{\phi}{\setU}{\setV}$ and $\funcAB{\psi}{\setU'}{\setV'}$ belong to $\Gsf$, then both $\funcAB{\phi^{-1}}{\setV}{\setU}$ and $\funcAB{\psi\circ\phi}{\phi^{-1}(\setU'\cap \setV)}{\psi(\setU'\cap \setV)}$ belong to $\Gsf$ (note that the diffeomorphism with empty domain, that is, the empty set, is allowed).
If $\Gsf$ is a set of locally defined smooth diffeomorphisms on $\M$, there exists a smallest local group of diffeomorphisms $\Gsf_\ast$ which contains $\Gsf$: we take the intersection $\cap\Fcal$ of the family $\Fcal$ of all local groups of diffeomorphisms which contain $\Gsf$ (note that $\Fcal$ is nonempty, since the set of all locally defined diffeomorphisms on $\M$ is such a local group). We call $\Gsf_\ast$ the \emph{local group of diffeomorphisms generated by $\Gsf$}.

Let $\Gsf$ be a local group of diffeomorphisms on $\M$. We define an equivalence relation on $\M$ by $x\sim y$ if $x=y$ or if there exists $\phi\in\Gsf$ such that $x\in\dom\phi$ and $\phi(x)=y$. The equivalence classes of this relation are called \emph{orbits of $\Gsf$}. Note that, if $x\in\M$ and there is no $\phi\in\Gsf$ such that $x\in\dom\phi$, the orbit of $x$ is the singleton of $x$. If $\Gsf$ is a set of locally defined smooth diffeomorphisms on $\M$, we define the \emph{orbits of $\Gsf$} as the orbits of $\Gsf_\ast$.

Given a locally defined smooth vector field $X$ on $\M$, we denote by $(X_t)_{t\in\Rn}$ the local one-parameter group of diffeomorphisms associated with $X$. If $\Dshf$ is a set of locally defined smooth vector fields on $\M$, we denote by $\Theta\Dshf$ the set of locally defined smooth diffeomorphisms on $\M$ given by
\begin{linenomath*}
	\begin{equation*}
\Theta\Dshf = \cup_{X\in\Dshf,t\in\Rn}X_t,
\end{equation*}
\end{linenomath*}
and by $\Psi\Dshf$ the local group of diffeomorphisms on $\M$ generated by $\Theta\Dshf$, i.e. the set of all finite compositions of local diffeomorphisms in $\Theta\Dshf$ (we are borrowing here the notation from \cite{Stefan},\cite{Stefan2}). We define the \emph{orbits of $\Dshf$} as the orbits of $\Theta\Dshf$. If $\vincl$ is a smooth distribution on $\M$, we define the \emph{orbits of $\vincl$} as the orbits of $\secloc(\vincl)$.   

We say that a smooth distribution $\vincl$ on $\M$ is \emph{invariant} by a local group of diffeomorphisms $\Gsf$ on $\M$ if, for each $x\in\M$, each $v\in\vincl_x$ and each $\phi\in\Gsf$ such that $x\in\dom\phi$, we have $\phi_\ast v\in\vincl_{\phi(x)}$, where $\phi_\ast$ denotes the tangent map of $\phi$.
We say that a smooth distribution $\vincl$ on $\M$ is \emph{invariant} by a set $\Gsf$ of locally defined smooth diffeomorphisms on $\M$ if it is invariant by $\Gsf_\ast$. We say that $\vincl$ is \emph{invariant} by a set $\Dshf$ of locally defined smooth vector fields on $\M$ if $\vincl$ is invariant by $\Psi\Dshf$. 

Given $\vincl$ and $\vincl'$ distributions on $\M$, we say that $\vincl\subset\vincl'$ if, for all $x\in\M$, $\vincl_x\subset\vincl_x'$. 

Given a smooth distribution $\vincl$ on $\M$ and a local group of diffeomorphisms $\Gsf$ on $\M$, there exists a smallest smooth distribution $\vincl^{\Gsf}$ on $\M$ which contains $\vincl$ and is invariant by $\Gsf$: if $\vincl$ is generated by the set of locally defined smooth vector fields $\Dshf$,  $\vincl^{\Gsf}$ is the distribution generated by the set of locally defined smooth vector fields  $\{\pushf{\phi}X\talque \phi\in\Gsf, X\in\Dshf\}$, where $\pushf{\phi}X$ denotes the pushforward of $X$ by $\phi$ (which is a locally defined smooth vector field on $\M$). Consequently, if $\Dshf$ is a set of locally defined smooth vector fields on $\M$, there exists a smallest smooth distribution $\Psf_\Dshf$ (this time we are borrowing the notation from \cite{Sussmann}) on $\M$ which contains the distribution $\lins{\Dshf}$ generated by $\Dshf$ and which is invariant by $\Dshf$, i.e. it is invariant by $\Psi\Dshf$. The smooth distribution $\Psf_\Dshf$ is generated by $\{\pushf{\phi}X\talque \phi\in\Psi\Dshf, X\in\Dshf\}$.

We can finally enunciate a version of the so-called \emph{orbit theorem}. The following statement is a subset of the the more general statements contained in \cite{Sussmann} (Theorem 4.1) and \cite{Stefan} (Theorems 1 and 5).

\begin{theorem}[orbit theorem]\label{thm:orbit}
	Let $\M$ be a finite dimensional paracompact smooth manifold and $\Dshf$ a set of locally defined smooth vector fields on $\M$. Then each orbit $S$ of $\Dshf$ is an immersed smooth submanifold of $\M$ such that, for each $x\in S$, the tangent space of $S$ at $x$ coincides with $\Psf_\Dshf(x)$. 
\end{theorem}

It was actually proved in \cite{Stefan} that each orbit $S$ of $\Dshf$ admits a unique smooth manifold structure which turns it into a \emph{leaf} of $\M$, i.e. a smooth immersed submanifold with the property that, for each locally connected topological space $N$ and each continuous map $\funcAB{f}{N}{\M}$ with image contained in $S$, the induced map $\funcAB{f}{N}{S}$ is continuous. Besides, the partition of $\M$ determined by the orbits of $S$ is a \emph{foliation with singularities} (cf. definition on page 700 of \cite{Stefan}). In particular, $\Psf_\Dshf$ is an involutive distribution (that was also proved in \cite{Sussmann}). It then follows that (recall that $\Dshf_\ast$ denotes the smallest involutive subset of locally defined smooth vector fields on $\M$ containing $\Dshf$) we have inclusions
\begin{linenomath*}
	\begin{equation*}
\lins{\Dshf}\subset\lins{\Dshf_\ast}\subset \Psf_\Dshf.
\end{equation*}
\end{linenomath*}
Indeed, the first inclusion is clear, and the second inclusion follows from the inclusion $\Dshf_\ast\subset\secloc(\Psf_\Dshf)$ (since, by the involutiveness of the distribution $\Psf_\Dshf$, $\secloc(\Psf_\Dshf)$ is an involutive set of locally defined smooth vector fields containing $\Dshf$, hence it must contain $\Dshf_\ast$) and from the fact that $\Psf_\Dshf$ is generated by $\secloc(\Psf_\Dshf)$. We therefore conclude that, if $\vincl$ is a smooth bracket-generating distribution on $\M$ and $\Dshf = \secloc(\vincl)$, then 
\begin{linenomath*}
	\begin{equation*}
\lins{\Dshf_\ast} = \Psf_\Dshf = \ft\M.
\end{equation*}
\end{linenomath*}
In particular, if $\M$ is connected, $\vincl$ admits a unique orbit which coincides with $\M$. We have thus proved the following version of Chow-Rashevskii's connectivity theorem. We say that a sectionally smooth curve on $\M$ is \emph{horizontal} with respect to $\vincl$ if all of its tangent vectors belong to $\vincl$.

\begin{corollary}[Chow-Rashevskii]\label{cor:Chow}
		Let $\M$ be a finite dimensional paracompact connected smooth manifold and $\vincl$ a smooth bracket-generating distribution on $\M$. Then $\M$ is \emph{$\vincl$-connected}, i.e. any two points in $\M$ may be connected by a sectionally smooth curve on $\M$ horizontal with respect to $\vincl$.
\end{corollary}

The converse to Chow-Rashevskii's theorem fails, i.e. the bracket-generating condition is not necessary for $\vincl$-connectivity (see \cite{Montgomery}, page 24). 

A necessary and sufficient condition for $\vincl$-connectivity may be obtained as a direct consequence of the following corollary of theorem \ref{thm:orbit} (cf. theorem 7.1 in \cite{Sussmann}).

\begin{corollary}[Sussmann's condition for $\Dshf$-connectivity]\label{cor:7.1Sussmann}
	Let $\M$ be a finite dimensional paracompact connected smooth manifold and $\Dshf$ a set of locally defined smooth vector fields on $\M$. Then $\M$ is \emph{$\Dshf$-connected} (i.e. $\M$ is an orbit of $\Dshf$) if, and only if,
	$$\Psf_\Dshf = \ft\M.$$
\end{corollary}

\subsection{Fiber and Parallel Derivatives}\label{subsec:Parallel}

Our last ingredient is a computational tool. Given a smooth linear subbundle $\vincl$ of $\ft\M$, we shall need to compute Lie brackets of vector fields in $\cv(\vincl)$. That could be accomplished by means of local charts on $\M$ and local trivializations of the vector bundle $\funcAB{\pi_\vincl}{\vincl}{\M}$, but in that case the computations we need to perform become rapidly messy. Instead, we compute by means of a method introduced in \cite{parallel} and summarized below.

Let $\funcAB{\pi_{E}}{E}{\M}$ be a smooth vector bundle over $\M$ endowed with a connection $\funcAB{\indconnection{}{}{E}}{\cv(\M)\times\secao(E)}{\secao(E)}$. The connection $\nabla^E$ defines a horizontal subbundle $\Hor(E)$ of $\ft E$, where $(\forall v_q\in E) \Hor_{v_q}(E)$ is the image of the \emph{horizontal lift at $v_q$}, $\funcAB{\lh_{v_q}}{\et{q}\M}{\et{v_q}E}$, defined by $w_q\mapsto \D V\cdot w_q$, where $\D$ denotes the tangent map and $V$ is any smooth local section of $\funcAB{\pi_E}{E}{\M}$ defined on an open neighborhood of $q$ such that $V(q)=v_q$ and $\indconnection{w_q}{V}{E}=0$. The horizontal lift $\funcAB{\lh_{v_q}}{\et{q}\M}{\et{v_q}E}$ is therefore a linear isomorphism onto $\Hor_{v_q}(E)$ whose inverse is the restriction of the tangent map $\D\pi_E$ to $\Hor_{v_q}(E)$. Denoting by $\Ver(E):=\ker\D\pi_E$ the \emph{vertical subbundle} of the tangent bundle $\ft E$, we thus have a Whitney sum decomposition
\begin{linenomath*}
	\begin{equation*}
\ft E = \Hor(E)\sd_E\Ver(E).
\end{equation*}
\end{linenomath*}
The \emph{connector} $\funcAB{\connector_E}{\ft E}{E}$ associated to the connection is given by $X_{v_q}\in\et{v_q}E\mapsto P_V(X_{v_q})\in\Ver_{v_q}(E)$ (where $P_V$ is the projection on the vertical subbundle induced by the Whitney sum decomposition above) followed by the inverse of the \emph{vertical lift} $\funcAB{\lv_{v_q}}{E_q}{\Ver_{v_q}(E)}$ at $v_q$ (which is the canonical linear isomorphism $E_q\equiv \et{v_q}(E_q)=\Ver_{v_q}(E)$). Note that, with these definitions:
\begin{enumerate}[1)]
	\item for all $X_{v_q}\in\ft E$, $X_{v_q} = \lh_{v_q}(\D\pi_E\cdot X_{v_q}) + \lv_{v_q}(\connector_E\cdot X_{v_q})$;
	\item for all $w_q\in\ft\M$, for all $V$ smooth local section of $\funcAB{\pi_E}{E}{\M}$ defined on an open neighborhood of $q$, we have $\indconnection{w_q}{V}{E}=\connector_E\cdot\D X\cdot w_q\in E_q$.
\end{enumerate}

Next, we consider two smooth vector bundles $\pi_E:E\rightarrow \M$ and $\pi_F:F\rightarrow \N$ 
over paracompact smooth manifolds $\M$ and $\N$, respectively, and 
$\funcAB{b}{E}{F}$ be a morphism of smooth fiber bundles (i.e. it
preserves fibers and is smooth) over $\funcAB{\tilde{b}}{\M}{\N}$. We denote by
$\funcAB{\fibdev b}{E}{\lin(E,\pullb{\tilde{b}}F)}$ the
\emph{fiber derivative} of $b$, i.e. the morphism of smooth fiber bundles defined by, for all
$v_q, w_q\in E_{q}$, $\fibdev b(v_q)\cdot w_q \doteq
\connector^V_F\cdot\D b\cdot \lv_{v_q}(w_q)\in F_{\tilde{b}(q)}$, where $\connector^V_F$ denotes the restriction of the connector $\kappa_F$ to the vertical subbundle (that is, $\connector^V_F$ is the inverse of the vertical lift). We don't need the connections to define the fiber derivative; what we need them for is to define the \emph{parallel derivative} $\pardev b:E\rightarrow\lin(\ft\M,\pullb{\tilde{b}}F)$. That is a smooth fiber bundle morphism defined by, for all $v_q\in E$ and
all $z_q\in\et{q}\M$,
\begin{linenomath*}
	\begin{equation*}
\pardev b(v_q)\cdot z_q \doteq  \connector_F\cdot\D
b\cdot\lh_{v_q}(z_q)\in F_{\tilde{b}(q)}.
\end{equation*}
\end{linenomath*}

The idea in considering these fiber and parallel derivatives is to provide a coordinate-free technique to
compute the tangent map of $b$, allowing its computation at a given element of $\ft E$ in terms of its vertical and horizontal
components, so that they play a role of ``partial derivatives''. That is to say, for all $X_{v_q}\in\ft E$, the following formulae hold:
\begin{linenomath*}
	\begin{equation*}
\begin{split}
\D\pi_F\cdot\D b\cdot X_{v_q} &= \D\tilde{b}\cdot\D\pi_E\cdot
X_{v_q}\\
\connector_F\cdot\D b\cdot X_{v_q} &= \fibdev
b(v_q)\cdot\connector_E\cdot X_{v_q} + \pardev
b(v_q)\cdot\D\pi_E\cdot X_{v_q}.
\end{split}
\end{equation*}
\end{linenomath*}

We finally come back to our initial setting, i.e. take $\M$ a finite dimensional paracompact smooth manifold endowed with a smooth linear subbundle $\vincl$ of $\ft\M$. We fix an auxiliary Riemannian metric tensor $\gtensor$ on $\M$, which induces a Whitney sum decomposition $\ft\M = \vincl\sd_\M\vincl^\perp$. We denote by $\funcAB{P}{\ft\M}{\vincl}$ the projection on the first factor determined by this Whitney sum, and by $\funcAB{\indconnection{}{}{\vincl}}{\cv(\M)\times\secao(\vincl)}{\secao(\vincl)}$ the connection on the vector bundle $\funcAB{\pi_\vincl}{\vincl}{\M}$ given by
\begin{linenomath*}
	\begin{equation*}
\indconnection{X}{Y}{\vincl} := P\connection{X}{Y},
\end{equation*}
\end{linenomath*}
where $\nabla$ is the Levi-Civita connection of $(\M,\gtensor)$.  Thus, both vector bundles $\funcAB{\tau_\M}{\ft\M}{\M}$ and $\funcAB{\pi_\vincl}{\vincl}{\M}$ are endowed with connections $\nabla$ (Levi-Civita) and $\nabla^\vincl$, with respective connectors and horizontal lifts denoted by $\connector, \lh_{v_q}$ and $\connector_\vincl, \lh_{v_q}^\vincl$. With respect to these connections, the Lie bracket $\lb{X}{Y}$ of (possibly locally defined) smooth vector fields, $X,Y\in\cv(\vincl)$ was computed in proposition 1 of \cite{parallel} by means of the following formulae, given $v_q\in\dom X\cap \dom Y$:
\begin{linenomath*}
	\begin{equation*}
\begin{split}
\connector_\vincl\cdot\lb{X}{Y}(v_q) &= \fibdev(\connector_\vincl\circ
Y)(v_q)\cdot\connector_\vincl\cdot X(v_q) + \pardev(\connector_\vincl\circ
Y)(v_q)\cdot\D\pi_\vincl\cdot X(v_q) -\\
&\quad -\fibdev(\connector_\vincl\circ X)(v_q)\cdot\connector_\vincl\cdot Y(v_q) -
\pardev(\connector_\vincl\circ X)(v_q)\cdot\D\pi_\vincl\cdot Y(v_q) +\\
&\quad +\ct^\vincl\bigpar{\D\pi_\vincl\cdot Y(v_q),\D\pi_\vincl\cdot X(v_q)}\cdot
v_q,\\
\D\pi_\vincl\cdot\lb{X}{Y}(v_q) &= \fibdev(\D\pi_\vincl\circ
Y)(v_q)\cdot\connector_\vincl\cdot X(v_q) + \pardev(\D\pi_\vincl\circ
Y)(v_q)\cdot\D\pi_\vincl\cdot X(v_q)-\\
&\quad - \fibdev(\D\pi_\vincl\circ
X)(v_q)\cdot\connector_\vincl\cdot Y(v_q) - \pardev(\D\pi_\vincl\circ
X)(v_q)\cdot\D\pi_\vincl\cdot Y(v_q),
\end{split}
\end{equation*}
\end{linenomath*}
where $\ct^\vincl$ is the curvature tensor of $\nabla^\vincl$.

We shall need the formulae above in the particular case in which: 1) $X$ is the nonholonomic vector field $\vf_{\vincl}$ of $(\M,\gtensor,\vincl)$, i.e. the vector field given by
\begin{linenomath*}
	\begin{equation*}
\vf_\vincl(v_q) = \lh_{v_q}^\vincl(v_q) = \D P\cdot \Sg(v_q),
\end{equation*}
\end{linenomath*}
where $\Sg$ is the geodesic spray of $(\M,\gtensor)$; 2) $Y$ is an arbitrary (locally defined) smooth vertical vector field. In this case, the above formulae simplify to, for all $v_q\in\dom Y$:
\begin{linenomath*}
	\begin{equation}\label{eq:lbhv}
\begin{split}
\connector_\vincl\cdot\lb{\vf_{\vincl}}{Y}(v_q) &= \pardev(\connector_\vincl\circ
Y)(v_q)\cdot v_q\\
\D\pi_\vincl\cdot\lb{\vf_\vincl}{Y}(v_q) &= -\connector_\vincl\cdot Y(v_q).
\end{split}
\end{equation}
\end{linenomath*}

\section{Statement and Proof of the Main Results}

\begin{theorem}[Smoothing in Chow's Theorem]\label{thm:ChowSmooth}
	Let $\M$ be a finite dimensional paracompact connected smooth manifold endowed with a smooth linear subbundle $\vincl$ of $\ft\M$. If $\vincl$ is bracket-generating, then any two points in $\M$ may be connected by a horizontal curve which is both a $\Ck{1}$ immersion and sectionally smooth, with arbitrary given initial and final velocities in $\vincl$.
\end{theorem}

\begin{proof}
	It suffices to consider the case $\dim\M\geq 2$, otherwise the thesis is trivial. Then, since $\vincl$ is bracket-generating, we must have $\rk\vincl\geq 2$; it then follows that the slit bundle $\vincl^\ast$ (i.e. $\vincl$ with the zero section removed) is a connected open submanifold of $\vincl$ (the fact that it is connected is a consequence of being the total space of a fiber bundle with fibers and base connected). We may apply the orbit theorem \ref{thm:orbit} to the paracompact connected smooth manifold $\vincl^\ast$ endowed with the set $\Dshf$ of locally defined smooth second order vector fields on $\vincl^\ast$, i.e. (noting that $\ft(\vincl^\ast) = \ft\vincl\rvert_{\vincl^\ast}$)
	\begin{linenomath*}
		\begin{equation*}
	\Dshf = \{X\in\secloc(\ft\vincl\rvert_{\vincl^\ast})\talque \forall v_q\in\dom X, \D\pi_\vincl\cdot X(v_q)=v_q\}.
	\end{equation*}
	\end{linenomath*}	
	We contend that $\Psf_\Dshf = \ft\vincl\rvert_{\vincl^\ast}$. Once we prove this contention, we conclude that each orbit of $\Dshf$ is a connected open submanifold of $\vincl^\ast$, which implies, due to the connectedness of $\vincl^\ast$, that $\vincl^\ast$ is the only orbit of $\Dshf$. That is to say, given $p,q\in\M$ and $v_p\in\vincl_p\setminus\{0\}$, $v_q\in\vincl_q\setminus\{0\}$, there exists a sectionally smooth curve in $\vincl^\ast$ connecting $v_p$ to $v_q$, whose smooth arcs are integral curves of vector fields in $\Dshf$, i.e. of second order vector fields. The projection on $\M$ of this sectionally smooth curve connects $p$ to $q$, with initial velocity $v_p$ and final velocity $v_q$, and it is both a sectionally smooth and a $\Ck{1}$-immersed horizontal curve on $\M$. By the arbitrariness of $p,q$ taken in $\M$ and of the initial and final velocities in $\vincl^\ast$, we have thus reached the thesis.
	
	It remains, therefore, to prove our contention, i.e. that $\Psf_\Dshf = \ft\vincl\rvert_{\vincl^\ast}$. Given $v_q\in\vincl^\ast$, we must prove that $\Psf_\Dshf(v_q)=\et{v_q}\vincl$, which will be done along the steps below. We fix an auxiliary Riemannian metric tensor $\gtensor$ on $\M$ and use the notation from subsection \ref{subsec:Parallel} of the preliminaries.
	
	\begin{enumerate}[1)]
		\item Since any local smooth vertical vector field in $\cv(\vincl^\ast)$ may be written as a difference of two smooth second order vector fields, i.e. of two vector fields in $\Dshf\subset\secloc(\Psf_\Dshf)$, and since $\Psf_\Dshf$ is a smooth distribution, we conclude that any local smooth vertical vector field in $\cv(\vincl^\ast)$ is a smooth local section of $\Psf_\Dshf$, which implies that the vertical space $\Ver_{v_q}(\vincl)$ is contained in $\Psf_\Dshf(v_q)$.
		
		\item Let $\vf_\vincl$ be the nonholonomic vector field of $(\M,\gtensor,\vincl)$ (which is a second order vector field in $\cv(\vincl)$, so that its restriction to the open submanifold $\vincl^\ast$ belongs to $\Dshf$) and $Y$ an arbitrary vertical smooth vector field in $\cv(\vincl^\ast)$ defined on an open neighborhood of $v_q$. Then both $\vf_\vincl\rvert_{\vincl^\ast}$ and $Y$ are sections of $\Psf_\Dshf$; since the latter smooth distribution is involutive, we conclude that the Lie bracket $\lb{\vf_\vincl}{Y}$ is a section of $\Psf_\Dshf$. But, as we have computed in \eqref{eq:lbhv}, $\D\pi_\vincl\cdot\lb{\vf_\vincl}{Y}(v_q) = -\connector_\vincl\cdot Y(v_q)$. It then follows that the vector
		\begin{linenomath*}
			\begin{equation*}
			\lh_{v_q}^\vincl\bigl(-\connector_\vincl\cdot Y(v_q)\bigr) = \lb{\vf_\vincl}{Y}(v_q) - \lv_{v_q}(\connector_\vincl\cdot\lb{\vf_\vincl}{Y}(v_q))
			\end{equation*}
		\end{linenomath*}
		belongs to $\Psf_\Dshf(v_q)$, as both vectors in the second member of the previous equality belong to that space. Since the restriction of $\connector_\vincl$ to $\Ver_{v_q}(\vincl)$ is a linear isomorphism onto $\vincl_q$ (it is the inverse of the vertical lift $\funcAB{\lv_{v_q}}{\vincl_q}{\Ver_{v_q}(\vincl)}$), and since the smooth vertical vector field $Y$ in $\cv(\vincl^\ast)$ on a neighborhood of $v_q$ was arbitrarily taken, we conclude that
		\begin{linenomath*}
			\begin{equation*}
			\lh_{v_q}^\vincl(\vincl_q)\subset \Psf_\Dshf(v_q).
			\end{equation*}
		\end{linenomath*}
	
		\item It follows from the previous step and from the arbitrariness of the fixed $v_q\in\vincl^\ast$ that, for any smooth locally defined vector field $X\in\secloc(\vincl)$, the horizontal lift  $X^\Hor\in\secloc(\ft\vincl\vert_{\vincl^\ast})$ defined by
		\begin{linenomath*}
			\begin{equation*}
			w_q \in\vincl^\ast\cap\pi_\vincl^{-1}(\dom X) \mapsto \lh_{w_q}^\vincl\bigl(X(q)\bigr)
			\end{equation*}
		\end{linenomath*}
		is a smooth local section of $\Psf_\Dshf$. Moreover, for all $w_q \in\vincl^\ast\cap\pi_\vincl^{-1}(\dom X)$, we have $\D\pi_\vincl\cdot X^\Hor(w_q) = X(q) = X\circ\pi_\vincl(w_q)$, i.e. the vector fields $X^\Hor$ and $X$ are $\pi_\vincl$-related. Then so are the Lie brackets of vector fields of this form, i.e. if $Y$ is another smooth locally defined vector field in $\secloc(\vincl)$, the locally defined vector fields $\lb{X^\Hor}{Y^\Hor}$ and $\lb{X}{Y}$ are $\pi_\vincl$-related. 
		
		As $\Psf_\Dshf$ is involutive, we conclude by induction on $k$ that, for an arbitrary $k$-tuple $X_1,\dotsc,X_k$ in $\secloc(\vincl)$ defined on an open neighborhood of $q$, $\lb{\dotsb [\lb{X_1^\Hor}{X_2^\Hor}}{\dotsb]X_{k-1}^\Hor], X_k^\Hor}$ is a smooth local section of $\Psf_\Dshf$ defined on a neighborhood of $v_q$ and the locally defined vector fields
		\begin{linenomath*}
			\begin{equation*}
			\lb{\dotsb [\lb{X_1^\Hor}{X_2^\Hor}}{\dotsb]X_{k-1}^\Hor], X_k^\Hor} \text{ and } \lb{\dotsb [\lb{X_1}{X_2}}{\dotsb]X_{k-1}], X_k}
			\end{equation*}
		\end{linenomath*}
		are $\pi_\vincl$-related. It then follows that the vector
		\begin{linenomath*}
			\begin{equation*}
			\begin{split}
			&\lh_{v_q}^\vincl\bigl( \lb{\dotsb [\lb{X_1}{X_2}}{\dotsb]X_{k-1}], X_k}(q)\bigr)=\\
			&\quad = \lb{\dotsb [\lb{X_1^\Hor}{X_2^\Hor}}{\dotsb]X_{k-1}^\Hor], X_k^\Hor}(v_q) -\\
			&\quad\quad -\lv_{v_q}\bigl(\connector_\vincl\cdot \lb{\dotsb [\lb{X_1^\Hor}{X_2^\Hor}}{\dotsb]X_{k-1}^\Hor], X_k^\Hor}(v_q)\bigr)
			\end{split}
			\end{equation*}
		\end{linenomath*}
		belongs to $\Psf_\Dshf(v_q)$, since both vectors on the second member of the previous equality belong to that space. 
		But, since $\vincl$ is a bracket-generating distribution, we have
		\begin{linenomath*}
			\begin{equation*}
			\et{q}\M = \linspan \{ \lb{\dotsb [\lb{X_1}{X_2}}{\dotsb]X_{k-1}], X_k}(q)\talque k\in\Nn,X_1,\dotsc,X_k\in \secloc(\vincl) \}.
			\end{equation*}
		\end{linenomath*}
		We finally conclude that $\Hor_{v_q}(\vincl) = \lh^\vincl_{v_q}(\et{q}\M)\subset\Psf_\Dshf(v_q)$. Thus, in view of step 1, we have
		\begin{linenomath*}
			\begin{equation*}
			\et{v_q}\vincl = \Hor_{v_q}(\vincl)\sd\Ver_{v_q}(\vincl)\subset\Psf_\Dshf(v_q),
			\end{equation*}
		\end{linenomath*}
		hence the equality holds in the above inclusion and our contention is proved.
	\end{enumerate}
\end{proof}

\begin{corollary}\label{cor:ChowSmooth}
	Let $\M$ be a finite dimensional paracompact smooth manifold endowed with a smooth linear subbundle $\vincl$ of $\ft\M$. If $\vincl$ is bracket-generating, then any two points belonging to a connected open subset $\setU\subset\M$ may be connected by a horizontal curve in $\setU$ which is both a $\Ck{1}$ immersion and sectionally smooth, with arbitrary given initial and final velocities in $\vincl$.
\end{corollary}

\begin{proof}
	Apply the previous theorem with $\setU$ in place of $\M$ and $\vincl\rvert_\setU$ in place of $\vincl$.
\end{proof}

We finally prove that the same smoothness property holds under Sussmann's condition for $\vincl$-connectivity (corollary \ref{cor:7.1Sussmann}).

\begin{theorem}[smoothness in Sussmann's condition for $\vincl$-connectivity]\label{thm:SussmannSmooth}
	Let $\M$ be a finite dimensional paracompact connected smooth manifold endowed with a smooth linear subbundle $\vincl$ of $\ft\M$ such that $\Psf_{\secloc(\vincl)} = \ft\M$. Then any two points in $\M$ may be connected by a horizontal curve which is both a $\Ck{1}$ immersion and sectionally smooth, with arbitrary given initial and final velocities in $\vincl$.
\end{theorem}

\begin{proof}
	As in the proof of theorem \ref{thm:ChowSmooth}, it suffices to consider the case $\dim\M\geq 2$, otherwise the thesis is trivial. Then, since $\Psf_{\secloc(\vincl)} = \ft\M$, we must have $\rk\vincl\geq 2$, so that the slit bundle $\vincl^\ast$ is a connected open submanifold of $\vincl$. Once more we consider the paracompact connected smooth manifold $\vincl^\ast$ endowed with the set $\Dshf$ of locally defined smooth second order vector fields on $\vincl^\ast$, i.e. 
\begin{linenomath*}
	\begin{equation*}
	\Dshf = \{X\in\secloc(\ft\vincl\rvert_{\vincl^\ast})\talque \forall v_q\in\dom X, \D\pi_\vincl\cdot X(v_q)=v_q\}.
	\end{equation*}
\end{linenomath*}	
We contend that $\Psf_\Dshf = \ft\vincl\rvert_{\vincl^\ast}$. Once we prove this contention, the thesis follows from Sussmann's condition \ref{cor:7.1Sussmann}.

Given $v_q\in\vincl^\ast$, we must prove that $\Psf_\Dshf(v_q)=\et{v_q}\vincl$, which will be done along the steps below.

\begin{enumerate}[1)]
	\item We fix an auxiliary Riemannian metric tensor $\gtensor$ on $\M$. Steps 1) and 2) in the proof of theorem \ref{thm:ChowSmooth} apply \emph{ipsis litteris}, so that both the vertical subspace $\Ver_{v_q}(\vincl)$ and the horizontal lift $\lh_{v_q}^\vincl(\vincl_q)$ are linear subspaces of $\Psf_\Dshf(v_q)$. Hence, for any smooth locally defined vector field $X\in\secloc(\vincl)$, the horizontal lift $X^\Hor\in\secloc(\ft\vincl\vert_{\vincl^\ast})$ is a smooth local section of $\Psf_\Dshf$.
	
	\item Since $\Psf_\Dshf$ is generated by $\secloc(\Psf_\Dshf)$, it follows from theorems 4.1 and 4.2 in \cite{Sussmann} that $\Psf_\Dshf$ is $\secloc(\Psf_\Dshf)$-invariant. Hence, for each $X\in\secloc(\vincl)$, we conclude from the previous step that $(X^\Hor_t)_{t\in\Rn}$ preserves $\Psf_\Dshf$.
	
	\item Let $w_q\in\et{q}\M$. Since $\et{q}\M = \Psf_{\secloc(\vincl)}(q)$, we may take $z_p\in\vincl$ and finite families $(X_i)_{1\leq i\leq k}$ of smooth local sections of $\vincl$ and $(t_i)_{1\leq i\leq k}$ of real numbers such that $(X_{k,t_k}\circ\dotsb\circ X_{1,t_1})_\ast z_p = w_q$. But, for any for any smooth locally defined vector field $X\in\secloc(\vincl)$, the horizontal lift  $X^\Hor\in\secloc(\ft\vincl\vert_{\vincl^\ast})$ is $\pi_\vincl$-related to $X$; it then follows, recalling that $\vf_\vincl$ denotes the nonholonomic vector field of $(\M,\gtensor,\vincl)$, that
	\begin{linenomath*}
		\begin{equation*}
		\begin{split}
		&\D\pi_\vincl\circ(X_{k,t_k}^\Hor\circ\dotsb\circ X_{1,t_1}^\Hor)_\ast\vf_\vincl(z_p) =\\
		&\quad = (X_{k,t_k}\circ\dotsb\circ X_{1,t_1})_\ast \circ\D\pi_\vincl\cdot \vf_\vincl(z_p)= w_q.
		\end{split} 
		\end{equation*}
	\end{linenomath*}

	We therefore conclude that
	\begin{linenomath*}
	\begin{equation*}
\begin{split}
	\lh_{v_q}^\vincl(w_q) &= (X_{k,t_k}^\Hor\circ\dotsb\circ X_{1,t_1}^\Hor)_\ast\vf_\vincl(z_p) -\\
	&\quad -\lv_{v_q}\bigl(\connector_\vincl\cdot(X_{k,t_k}^\Hor\circ\dotsb\circ X_{1,t_1}^\Hor)_\ast\vf_\vincl(z_p)\bigr).
\end{split}
	\end{equation*}
\end{linenomath*}
Hence, $\lh_{v_q}^\vincl(w_q)$ belongs to $\Psf_\Dshf(v_q)$, since both vectors on the second member of the previous equality belong to that space, in view of steps 1 and 2.
Since $w_q\in\et{q}\M$ was arbitrarily taken, we conclude that $\Hor_{v_q}(\vincl) = \lh^\vincl_{v_q}(\et{q}\M)\subset\Psf_\Dshf(v_q)$. Thus, $\et{v_q}\vincl = \Hor_{v_q}(\vincl)\sd\Ver_{v_q}(\vincl)\subset\Psf_\Dshf(v_q)$, hence the equality holds in the above inclusion and our contention is proved.
	
\end{enumerate}

\end{proof}

\bigskip
%\begin{center}{\scshape  Acknowledgements} \end{center}
%\clearpage

\bibliographystyle{siam}
%\addcontentsline{toc}{chapter}{Bibliografia}
\bibliography{biblio}

% assinatura
%\cleardpage \thispagestyle{empty} \vspace*{2cm}

%\vspace{1cm}\begin{ass3}{\asslen}{\rulelen}{\minilen}
%\eu\\
%\vspace*{1cm}%\rule{\tamrule}{0.2 mm}\\
%\reinaldo\\
%\vspace*{1cm} São Paulo, 18 de julho de 2017.
%\end{ass3}

\end{document}